\documentclass[11pt,a4paper]{article}

\usepackage{amsmath,amssymb,amsthm,graphics,amscd}
\usepackage{longtable}
\usepackage{cite}
\usepackage{color}
\usepackage{graphicx}
\usepackage{epsf}
\usepackage{fancyhdr}
\usepackage{pxfonts}

\long\def\symbolfootnote[#1]#2{\begingroup\def\thefootnote{\fnsymbol{footnote}}
\footnote[#1]{#2}\endgroup}
  \renewenvironment{thebibliography}[1]{
    \begin{oldthebibliography}{#1}
      \setlength{\parskip}{0ex}
      \setlength{\itemsep}{0ex}
  }
  {
    \end{oldthebibliography}
  }
  
\setlongtables
\begin{document}

\newcounter{rownum}
\setcounter{rownum}{0}
\newcommand{\ab}{\addtocounter{rownum}{1}\arabic{rownum}}
\newcommand{\im}{\mathrm{im}}
\newcommand{\x}{$\times$}

\newtheorem{lemma}{Lemma}[section]
\newtheorem{theorem}[lemma]{Theorem}
\newtheorem*{ttt}{Theorem}
\newtheorem*{T1}{Theorem 1}
\newtheorem*{T2}{Theorem 2}
\newtheorem{corollary}[lemma]{Corollary}
\newtheorem{conjecture}[lemma]{Conjecture}
\newtheorem{prop}[lemma]{Proposition}
\newtheorem*{question}{Question}
\newtheorem*{questions}{Questions}
\newtheorem{q}[lemma]{Question}
\theoremstyle{remark}
\newtheorem{remark}[lemma]{Remark}
\theoremstyle{definition}
\newtheorem{defn}[lemma]{Definition}
\newtheorem{example}[lemma]{Example}
\newtheorem*{claim}{Claim}

\renewcommand{\labelenumi}{(\roman{enumi})}
\newcommand{\Hom}{\mathrm{Hom}}
\newcommand{\Ext}{\mathrm{Ext}}
\newcommand{\soc}{\mathrm{Soc}}

\newenvironment{changemargin}[1]{%
  \begin{list}{}{%
    \setlength{\topsep}{0pt}%
    \setlength{\topmargin}{#1}%
    \setlength{\listparindent}{\parindent}%
    \setlength{\itemindent}{\parindent}%
    \setlength{\parsep}{\parskip}%
  }%
  \item[]}{\end{list}}
\addtolength{\parskip}{0.5\baselineskip}
\parindent=0pt

\title{Unbounding Ext}
\author{David I. Stewart}
\date{New College, Oxford}
\maketitle

{\bf Abstract.} We produce examples in the cohomology of algebraic groups which answer two questions of Parshall and Scott. Specifically, if $G=SL_2$, then we show: (a) $\dim \Ext_G^2(L,L)$ can be arbitrarily large for a simple module $L$; and (b) if we define $\gamma_m=\max_L\dim H^m(G,L)$ where the maximum is taken over all simple $G$-modules $L$, then the sequence $\{\gamma_m\}$ grows exponentially fast with $m$.

\section*{Introduction}
Let $G$ be a simply connected, semisimple algebraic group with associated root 
system $\Phi$ defined over an algebraically closed field $k$ of characteristic $p>0$. We mention some notation taken to be consistent with \cite{Jan03}; any undefined notation can be found in there. Let $B$ be a Borel subgroup of $G$ with maximal torus $T$ defining a set of dominant weights $X^+(T)$, a subset of the weight lattice $X(T)$ of $T$, where $X(T)\cong \mathbb Z^n$; if $\lambda\in X(T)$ we write $\lambda=(a_1,a_2\dots,a_n)$. Recall that the simple $G$-modules are indexed by highest weight $\lambda\in X^+(T)\cong \mathbb Z_{\geq 0}^n$, the modules are then denoted by $L(\lambda)$. In the case $G=SL_2$ we identify $X^+(T)$ with $\mathbb Z_{\geq 0}$. Let $X_1(T)$ denote the $p$-restricted weights; that is the set of $(a_1,\dots,a_n)=\lambda\in X(T)$ with each $a_i<p$. Then any weight $\lambda\in X(T)$ has a $p$-adic expansion $\lambda=\lambda_0+p\lambda_1+\dots+p^n\lambda_n$, for some $n\in \mathbb N$ with each $\lambda_i\in X_1(T)$. We denote by $X_{e,p}$ the subset of $X^+(T)$ consisting of weights whose $p$-adic expansion is no longer than $e$, that is $X_{e,p}=\{\lambda\in X^+(T):\lambda_r=0 ,\ \forall r>e\}$.

In \cite{PS11} the authors find a constant $c:=c(\Phi, n, e)$ such that 
\[\dim \Ext^n_G(L(\lambda), L(\mu))\leq c\] 
for all simply connected, semisimple algebraic groups with root system $\Phi$ (thus 
independent of the characteristic $p$ of $k$) and all $\lambda\in X_{e,p}$. 

In the case $n = 1$, the authors are able to drop the dependence on $e$ to yield a 
constant $c:=c(\Phi)$ such that 
\[\dim \Ext^1_G(L(\lambda), L(\mu))\leq c\]
for all simply connected, semisimple algebraic groups with root system $\Phi$. In 
[ibid., Remark 7.4(b)] the authors ask if the dependence on the length $e$ of 
the $p$-adic expansion of $\lambda$ can be dropped for $n > 1$. 

Let $p > 2$ and let $G = SL_2$. In Theorem 1 we give a sequence of weights $\lambda_r,\mu_r\in X_+(T)$
for $G$ such that $\dim \Ext^2_G(L(\lambda_r), L(\mu_r)) = r$, answering this question in the negative. This is the subject of Section 1.

In a further paper, \cite{PS10}, the same authors make the following definitions: For 
an algebraic group $G$ and (rational) $G$-module $V$, put
\begin{align*}\gamma_m(V)&=\max_{L-\text{irred}}\dim\Ext^m_G(V,L)\\
\gamma_m(\Phi,e,p)&=\max_{\lambda\in X_{e,p}}\gamma_m(L(\lambda))\\
\gamma_m(\Phi,e)&=\max_p\gamma_m(\Phi,e,p),\end{align*}
where the maximum in the first line is over all irreducible $G$-modules $L$. These are finite by \cite[7.1]{PS11}. They prove

\begin{theorem}[{\cite[6.1]{PS10}}]\label{PS} \begin{enumerate}\item The sequence $\{\log \gamma_m(\Phi,e)\}$ has polynomial rate of growth at most $4$.
\item For any fixed prime $p$, the sequence 
$\{\log \gamma_m(\Phi, e, p)\}$ has polynomial rate of 
growth at most 3. \end{enumerate}\end{theorem}

They then ask if these bounds can be improved to polynomial rates of growth in the case of cohomology. To wit, the following is Question 6.2 in [ibid.]: 

\begin{q}Let $\Phi$ be a finite root system. Do there exist constants $C=C(\Phi)$ 
and $f = f(\Phi)$ such that \[\dim H^m(G, L) \leq Cm^f\] 
for all semisimple, simply connected groups $G$ over an algebraically closed field 
$k$ (of arbitrary characteristic) having root system $\Phi$ and all irreducible rational 
$G$-modules $L$? \end{q}

Let again $G = SL_2$ and let $p$ be arbitrary. Define $\gamma_m=\max_{L-\text{irred}}\dim H^m(G, L)$, again with the maximum over all irreducible $G$-modules $L$. We use the algorithm in \cite{Par07} to show that the sequence $\{\gamma_m\}$ grows exponentially with $m$, answering this second question in the negative. For simplicity we prove this first in the case $p=2$. Recall that there is a Frobenius map $F:G\to G$; induced by raising matrix entries to the $p$th power. Composing $F$ with a representation $G\to GL(V)$ gives a new $G$-module $V^{[1]}$ whose weights are $p$ times the weights of $V$. We show that the sequence $H^m(G,L(1)^{[m]})=\Pi_{m-1}$ where $L(1)^{[m]}$ is the $m$th Frobenius twist of the natural module $L(1)$ for $G$ and $\Pi_{m}$ is the number of partitions of unity into $m$ powers of $1/2$.  This is our Theorem 2. We prove this in Section 2 and offer a number of extensions to this result, including to the case $p>2$. 

In D. Hemmer's MathSciNet review of \cite{Par07}, he admits to being unsure how difficult the recursions would 
be to use for actual computation. We hope our theorem serves as a vindication of the usefulness of Parker's algorithm for producing interesting general 
results about the behaviour of $\Ext$-groups. 

At the end of the paper, we make a number of remarks indicating, as far as we can, various possible extensions to this work. We also make some remarks of relevance to questions of \cite{GKKL} which considers the putative existence of bounds on the dimension of the cohomology group $H^n(G,V)$ in terms of (powers of) the dimension of $V$, where $G$ is a finite group and $V$ an absolutely irreducible $kG$-module.

{\bf Acknowledgements} We thank Len Scott for helpful comments and suggestions while this paper was being produced. We also thank Chris Bowman for some helpful telephone conversations. Our thanks also to the referee for suggesting improvements.

\section{Unbounding $\Ext$}

Let $G=SL_2$ defined over an algebraically closed field k of characteristic $p>3$. 
The following result is the main result from \cite{Ste10}. 

\begin{lemma}Let $V=L(r)^{[d]}$ be any Frobenius twist (possibly trivial) of the 
irreducible $G$-module $L(r)$ with highest weight $r$ where $r$ is one of 
\begin{align*}&2p \\
&2p^2-2p-2\\
&2p-2 + (2p-2)p^e \ (e > 1) \end{align*}
Then $H^2(G, V ) \cong k$. For all other irreducible $G$-modules $V$, $H^2(G, V) = 0$. 
\end{lemma}

Now we can prove 
\begin{T1}Let $V_n = L(1)\otimes L(1)^{[1]}\otimes\dots\otimes L(1)^{[n]}$. 

Then $\Ext^2_G(V_n,V_n)=n.$ \end{T1}
\begin{proof}By Steinberg's tensor product theorem, $V_n$ is simple; thus it is self-dual 
and we have \begin{align*}{\small \Ext^2_G(V_n , V_n )&\cong\Ext^2_G(k, V_n\otimes V_n^*)\\ 
&\cong H^2 (G, (L(1)\otimes L(1))\otimes (L(1)\otimes L(1))^{[1]}\otimes\dots\otimes (L(1)\otimes L(1))^{[n]}) \\
&\cong H^2(G, (L(2) \oplus k)\otimes  (L(2) \oplus k)^{[1]}\otimes\dots\otimes (L(2) \oplus k)^{[n]}) \\
&\cong H^2(G, k) \oplus H^2(G, L(2)) \oplus H^2 (G, L(2)^{[1]})\oplus\dots\oplus H^2(G, L(2)^{[n]})
\\&\oplus H^2 (G, L(2)\otimes  L(2)^{[1]} ) \oplus\dots\ .}\end{align*}
The third isomorphism follows since when $p>2$, $L(1)\otimes L(1)$ has composition factors $L(2)$ and $k$ which do not extend each other. The last isomorphism is a formal expansion of the tensor product in the third line, using the fact that the Frobenius twist, tensor product and the functors $H^i(G,?)$ commute with direct sums; the modules $L(2)^{[i_1]}\otimes L(2)^{[i_2]}\otimes \dots\otimes L(2)^{[i_r]}$ for distinct $i_j$ are simple by Steinberg's tensor product theorem.

Now, by the Lemma, the only terms in this expression which are non-zero are $H^2 (G, L(2)^{[d]})$ 
with $d > 0$. Thus $\dim\Ext^2_G(V_n , V_n ) = n$ as required. \end{proof}

\begin{remark}In fact one knows from \cite{McN02} that if $p\geq h$, then for any $r > 0$, we have $H^2(G, \mathfrak{g}^{[r]}) \cong k$ for any simply connected, simple algebraic group $G$, 
where $\mathfrak g$ denotes the Lie algebra of $G$.

Then one can construct a similar example to the above for any $G$. One takes any simple module $L=L(\lambda)$ such that $L$ is a faithful representation of $G$, with $p$ big enough so that $L 
\otimes  L^*$ is completely reducible. Then it will contain $\mathfrak g$ and $k$ as direct summands. If the weights of $L 
\otimes  L^*$ are 
less than $p^r$ then one can take $V_n = L\otimes L^{[r]}\otimes L^{[2r]}\dots L^{[nr]}$ with the property 
that $\dim\Ext^2(V_n , V_n )\geq n$.
\end{remark}

We now know that Parshall and Scott's restriction on the length of the $p$-adic expansion of $L$ is necessary to have a finite bound for $\max\dim \Ext_G^n(L,L')\leq c(\Phi,e)$ with the maximum taken over all irreducible modules $L,L'$ with $e_p(L)<e$. In which case, it might be interesting to see how the sequence \[\{f_{e}\}:=\max\{\dim \Ext_G^n(L,L')\}\] grows with $e$ for fixed values of $n$ and $\Phi$, where the maximum is taken over all $p$ and irreducible $G$-modules $L,\ L'$ with $e_p(L)<e$. In the case $n=2$ our examples show that $f_e$ is at least linear.

\section{Exponential growth of $H^m$}

Let $G = SL_2$ defined over an algebraically closed field $k$ whose characteristic will be $p = 2$ until further notice (i.e. Remark \ref{changeg}). 

In this section we show that the sequence $\{\dim H^n(G, L(2^n))\}$ has exponential growth with $n$. (In fact, it is true that $\dim H^n(G, L(2^n))=\max_m \dim H^n(G,L(2^m))$, see Remark \ref{generic} below.)

We need the following two formulae from \cite{Par07}, valid when $p = 2$. 

\begin{theorem}\label{parkereq} Let $M$ be a $G$-module and take $b, q\in \mathbb N$ with $q > 0$. Then 
\begin{align}
\Ext^q_G(\Delta(2b),M^{[1]})&\cong\bigoplus_{n=0}^{n=q}\Ext^{q-n}_G(\Delta(n + b), M),\\
\Ext^q_G(\Delta (2b + 1), M^{[1]}\otimes L(1))&\cong\Ext^q_G(\Delta(b),M),\end{align}\end{theorem}
where $\Delta(r)$ denotes the Weyl module for $G$ of highest weight $r$.

Note that the above formulae are also clearly valid when $q = 0$; however, our 
analysis of the algorithm is slightly more transparent if we do not use these 
formulae in the case $q = 0$. 

Using (1) and (2) it is possible to calculate $H^q (G, L)$ inductively for any simple 
$G$-module $L$. We give such a recipe now. 

Firstly, by Steinberg's Tensor Product Theorem, $L\cong L(a_0 )\otimes L(a_1)^{[1] }\otimes L(a_2 )^{[2]}\otimes \dots \otimes L(a_n)^{[n]}$ for some $n\in \mathbb N$, with (as $p=2$) each $a_i\in\{0, 1\}$, i.e. $L$ is the trivial module $k$, or 
a tensor products of different Frobenius twists of the natural module $L(1)$ for $G$. By the 
linkage principle, if $H^q(G, L)\neq 0$, then $L = M^{[1]}$ for some simple module $M$. 

Thus, taking $b = 0$, we apply (1) to express $H^q(G,M^{[1]})\cong \Ext^q_G(k,M^{[1]})$ in terms of $\Ext$s of equal or 
lower degree between $\Delta$-modules and another simple module $M$ of lower weight.

We may then ignore about half of these $\Ext$ terms since, if the parities of the 
highest weights of $M$ and a given $\Delta (r)$ module are different then this $\Ext$ 
term vanishes by linkage. For the remainder, apply equation (2) if $M$ is a 
simple module of odd high weight; and then continue to expand each surviving 
$\Ext$ term using equation (1). Eventually this process terminates with a sum 
of terms $\Ext^q_G(\Delta (r), k)$ with $q > 0$, which are $0$ by \cite[II.4.13]{Jan03} and terms $\Ext^0_G(\Delta (r_i ), N_i ) \cong \Hom_G(\Delta (r_i ), N_i )$ for some known collection of simple modules $N_i$. We call these $\Ext^0$ terms \emph{leaves}; see below for an example.

As each $N_i$ is simple and $\Delta(r_i )$ has a simple head, each of these leaves is then visibly either isomorphic to $k$ or $0$ (according to whether or not the highest weight of $N_i$ is the integer $r_i$) and so the desired value of $\dim H^q(G, L)$ has been calculated.

Given a simple module $L$ and a degree $m$ of cohomology, we wish to enumerate these $\Ext^0$ leaves. To this end we make the following  recursive definition, which will be elucidated by the following examples.

\begin{defn}\label{def}
For a given degree $m>0$ of cohomology, and simple module $L$, define an $a$-string to be a list of non-negative integers $(a_1,\dots, a_n)$ with $a_n>0$ and $\sum a_i=m$ such that the following procedure terminates successfully.

Set $\mathcal T_1$ to be the term $E=\Ext_G^m(\Delta(0),L)$. 

At stage $1$, if the parity of the highest weight of $L$ is odd, then return failure. Otherwise use equation (1) to expand $\mathcal T_1=E$, and then consider the term $\Ext_G^{m-a_1}(\Delta(a_1),L^{[-1]})$. If $n=1$, (and so $m=a_1$) then terminate, returning the leaf `$\Ext_G^{0}(\Delta(a_1),L^{[-1]})$'. Otherwise set $\mathcal T_2=\Ext_G^{m-a_1}(\Delta(a_1),L^{[-1]})$ and continue to step $2$. 

At stage $r$, one is given $\mathcal T_r=\Ext_G^{m-\sum_{i=1}^{r-1}(a_i)}(\Delta(x),L(y))$ for some $x,y\in \mathbb N$. Check the parities of $x$ and $y$. If they are different then return failure; otherwise, if necessary, apply (2) to $\mathcal T_r$ and replace it with the resulting term, until either the parities of the weights in $\mathcal T_r$ differ, whence return failure, or until they are both even. Then use equation (1) to expand $\mathcal T_r$ and consider the resulting term $\Ext_G^{m-\sum_{i=1}^{r}a_i}(\Delta(x'),L(y'))$  for some $x',y'\in \mathbb N$. If $r=n$ then terminate, returning the leaf  $\Ext_G^{0}(\Delta(x'),L(y'))$. Otherwise set $\mathcal T_{r+1}=\Ext_G^{m-\sum_{i=1}^{r+1}a_i}(\Delta(x'),L(y'))$ and continue to step $r+1$.
\end{defn}

\begin{example}Let $m=6$ and $L=L(24)$. Then there is an $a$-string $(4,0,2)$:
\begin{align*}{\small \Ext^6(\Delta(0),L(24))&\cong\Ext^6(\Delta(0),L(12))\oplus \Ext^5(\Delta(1),L(12))\oplus \Ext^4(\Delta(2),L(12))\\&\oplus \Ext^3(\Delta(3),L(12))\oplus \underline{\Ext^2(\Delta(4),L(12))}\oplus \Ext^1(\Delta(4),L(12))\\&\oplus \Ext^0(\Delta(6),L(12))\\
\Ext^2(\Delta(4),L(12))&\cong\underline{\Ext^2(\Delta(2),L(6))}\oplus \Ext^1(\Delta(3),L(6))\oplus \Ext^0(\Delta(4),L(6))\\
\Ext^2(\Delta(2),L(6))&\cong\Ext^2(\Delta(1),L(3))\oplus \Ext^1(\Delta(2),L(3))\oplus \underline{\Ext^0(\Delta(3),L(3))}\\
\Ext^0(\Delta(3),L(3))&\cong k, }\end{align*}where we have underlined the terms corresponding to the $a_i$. In this case the $a$-string happens to give a non-trivial leaf, showing in particular, that $\Ext_G^6(L(0),L(24))>0$.\end{example}

Note that not all strings of non-negative integers adding up to $m$ are valid $a$-strings. For instance, in the setting of the above example, strings such as $(3,3)$ or $(3,2,1)$ are not $a$-strings since they gives rise to a chain
\[\dim\Ext^6(\Delta(0),L(24))\geq \dim\Ext^3(\Delta(3),L(12))=0,\]
as the parity of $3$ and $12$ is different so the procedure of the definition returns failure.

Also there are $a$-strings which result in $\Ext^0$ leaves which are zero. For instance the string $42$ is valid as an $a$-string:

\[\dim \Ext^6(\Delta(0),L(24))\geq \dim \Ext^2(\Delta(4),L(12))\geq \dim\Ext^0(\Delta(4),L(6))\]

but zero. We call an $a$-string which results in a non-zero leaf, a non-trivial $a$-string. Thus we have $\dim H^m(G,L)=|\{\text{non-trivial $a$-strings}\}|$. We wish to give a lower bound on the number of non-trivial $a$-strings. 

Firstly though, define an $(a,n)$-string to be a string of length $n$ so that the first $r$ entries are an $a$-string (of length $r\leq n$) and the remaining entries are $0$. We can of course, recover the original $a$-string from an $(a,n)$-string by removing all $0$s from the end. If the highest weight of $L$ is no more than $2^n$ then the length of any valid a-string can be no longer than $n$ and so we have a bijection between $a$-strings and $(a,n)$-strings.

Let $L=L(2^n)$. It is clear that procedure can be applied a maximum of $n$ times. So all $a$-strings for this $L$ can be made into $(a,n)$-strings. Keeping $L=L(2^n)$, we have the
\begin{lemma}An $(a,n)$-string $(a_1,\dots, a_n)$ is non-trivial provided there exists a string of positive integers $(b_1,b_2,\dots,b_n)$ with \begin{enumerate}\item $a_i+b_{i-1}=2b_i$ for $1\leq i\leq n-1$, \item $b_n=a_n$ and \item $b_{n-1}+b_n=1$,\end{enumerate}\end{lemma}
where we also set $a_0=b_{-1}=b_0=0$.\begin{proof}To prove the lemma, we trace the highest weight of the $\Delta$-module on the left through the procedure given in Definition \ref{def}. One finds that at step $r<n$ one is given a term 
\[\mathcal T_r=\Ext_G^{m-\sum_{i=1}^{r-1} a_i}\left(\Delta\left(\frac{\frac{\frac{\frac{a_1}{2}+a_2}{2}+\dots}{\ddots}\ddots+a_{r-2}}{2}+a_{r-1}\right),L(2^{n-r+1})\right).\]
Since the $(a,n)$-string is assumed to be non-trivial the parity of the left hand side must be even. Inductively, assume that we have defined an integer \[b_{r-2}=\frac{\frac{\frac{\frac{a_1}{2}+a_2}{2}+\dots}{\ddots}\ddots+a_{r-2}}{2}.\] Since the highest weight of the $\Delta$-module, $a_{r-1}+b_{r-2}$, is even, we may set $a_{r-1}+b_{r-2}=2b_{r-1}$. Thus \[b_{r-1}=\frac{\frac{\frac{\frac{a_1}{2}+a_2}{2}+\dots}{\ddots}\ddots+a_{r-1}}{2}\] as required. Finally, taking $r=n$ and expanding one last time we have a term $\Ext^0(\Delta(b_{n-1}+b_n),L(1))$ which is non-zero (and one-dimensional) precisely if $b_{n-1}+b_n=1$ as required.\end{proof}

As each $a_i$ is positive, the resulting string has the property (i'): $2b_i\geq b_{i-1}$ for $2\leq i\leq n-1$. Note also that if such a string exists for a given non-trivial $(a,n)$-string, it has property (iv): $m=\sum a_i=(\sum_{i=1}^n b_i)+b_{n-1}$; so $\sum_{i=1}^{n-1}b_i=m-1$. We call a string satisfying properties (i'), (iii) and (iv) a $b$-string, and observe that if a $b$-string exists for a given $(a,n)$-string, one can recover the original $a$-string.

Indeed, the proof of the lemma shows that any $b$-string gives rise to a non-trivial $a$-string. So it suffices to count $b$-strings. We do this now in the case $n=m-1$.

Take $n=m-1$. If $b_{n-1}=0$ then $b_1=\dots=b_{n-2}=0$ by property (i'); thus $m=1$ by property (iv) and thus $n=m-1=0$ which is nonsense. So $b_{n-1}=1$. Then we wish to find all sequences $b_1,\dots,b_{n-2}$ with $\sum b_i=n-2$ and $b_i\geq 2b_{i-1}$. Reversing the order; call a string of $n-1$ integers a $(c,n-1)$-string if $c_1=1$ and $c_i\leq 2c_{i-1}$ with $\sum_{i=1}^{n-1} c_i=n-1$. For each $n$, set $\Pi_{n-1}$ equal to the number of $(c,n-1)$-strings; this is then  precisely the sequence $H_{n-1}$ from \cite[p150]{FP87}. Thus we have that the dimension of $H^m(G,L(2^m))$ is the integer $\Pi_{m-1}$: the number of `level number sequences' associated to binary trees, or the number of partitions of $1$ into $m$ powers of $1/2$.\footnote{See \cite[http://oeis.org/A002572]{OEIS} for more on this sequence} We have from \cite{FP87} the inequality \[F_n\leq H_n\leq 2^{n-1}.\]

As $F_n\sim \left(\frac{1+\sqrt(5)}{2}\right)^n$, it follows immediately that $H_n$ grows exponentially, but we give a quick proof here that $\Pi_{2n+1}\geq 2^n$: 

Observe
\[1,\underbrace{2,2,\dots,2}_{n}\underbrace{0,0,\dots,0}_{n}\]
is a $c$-string. For any choice of subset of the $2$s in the first underbrace, we may replace each $2$ by the string $1,1$ and remove a $0$ from the right to have another $c$-string. Running through the different choices of the $2^{n}$ subsets we see that they are all distinct; and thus 

\begin{T2}\label{mine}For $m>2$, \[\dim H^{2m}(G,L(2^{2m}))\geq 2^{m-1}\] and so $H^{m}(G,L(2^{m}))$ grows exponentially with $m$.\end{T2}

\begin{remark}\label{generic} The longest $b$-string without $0$s at the front is clearly \[\underbrace{1,1,\dots,1}_{m-1},0.\] It follows from this that $\dim H^m(G,L(1)^{[r]})<\dim H^m(G,L(1)^{[m]})$ if and only if $r<m$, with equality otherwise. So for $p=2$, rational stability occurs for the module $L(1)$ at the Frobenius twist $m$, in other words the value of $\epsilon$ in \cite[Corollary 6.8]{CPSK} can be as large as $m$.

Since the dimensions of rationally stable and generic cohomology $H_\mathrm{gen}$ are a common limit, this shows in particular that when $p=2$, we have \[\dim H^m_{\mathrm{gen}}(G,L(1))=\Pi_{m-1}.\]\end{remark}

\begin{remark} We note that the rate of growth of $H^m$ is not too severely underestimated by a sequence $\{C.2^{m/2}\}$. The following are the precise numbers up to $n=31$:

\begin{figure}[ht]
\begin{minipage}[b]{0.5\linewidth}{\footnotesize
\centering

\begin{verbatim}H^4(G,L(2^4))=2
H^5(G,L(2^5))=3
H^6(G,L(2^6))=5
H^7(G,L(2^7))=9
H^8(G,L(2^8))=16
H^9(G,L(2^9))=28
H^10(G,L(2^10))=50
H^11(G,L(2^11))=89
H^12(G,L(2^12))=159
H^13(G,L(2^13))=285
H^14(G,L(2^14))=510
H^15(G,L(2^15))=914
H^16(G,L(2^16))=1639
H^17(G,L(2^17))=2938\end{verbatim}}
\end{minipage}
\hspace{0.5cm}
\begin{minipage}[b]{0.5\linewidth}{\footnotesize
\centering
\begin{verbatim}H^18(G,L(2^18))=5269
H^19(G,L(2^19))=9451
H^20(G,L(2^20))=16952
H^21(G,L(2^21))=30410
H^22(G,L(2^22))=54555
H^23(G,L(2^23))=97871
H^24(G,L(2^24))=175586
H^25(G,L(2^25))=315016
H^26(G,L(2^26))=565168
H^27(G,L(2^27))=1013976
H^28(G,L(2^28))=1819198
H^29(G,L(2^29))=3263875
H^30(G,L(2^30))=5855833
H^31(G,L(2^31))=10506175\end{verbatim}}
\end{minipage}
\end{figure}
Indeed \cite[Theorem 1]{FP87} shows that $H_\text{gen}^m(G,L(1))=\Pi_{m-1}\sim K\nu^m$, where $K\sim 0.255$ and $\nu\sim 1.794$ are constants defined in [{\it ibid.}] 

This would suggest that the best  likely result in the spirit of Theorem \ref{PS} given in the introduction would be that the sequence $\{\log \gamma_m(\Phi,e)\}$ has polynomial growth at most $1$ for any $\Phi$ (in other words, is linear with $m$). In any case, Theorem 2 shows that Parshall and Scott's estimate is certainly in the right ball-park.
\end{remark}

\begin{remark}\label{changeweight} One can replace the weight $2^{m}$ with any other weight $r.2^{m}$ with the result that the sequence $\{\dim H^m(SL_2,L(r.2^{m}))\}$ grows exponentially fast. We have written a computer program using Parker's algorithm to calculate the dimensions of cohomology groups. The output from the program giving dimensions for $H^m(SL_2,L(r.2^{m-2}))$ is given below.

\begin{figure}[ht]
\begin{minipage}[b]{0.5\linewidth}{\footnotesize
\centering

\begin{verbatim}H^3(G,L(3.2))=1
H^4(G,L(3.2^2))=1
H^5(G,L(3.2^3))=2
H^6(G,L(3.2^4))=4
H^7(G,L(3.2^5))=6
H^8(G,L(3.2^6))=11
H^9(G,L(3.2^7))=20
H^10(G,L(3.2^8))=35
H^11(G,L(3.2^9))=63
H^12(G,L(3.2^10))=113
H^13(G,L(3.2^11))=201
H^14(G,L(3.2^12))=361
H^15(G,L(3.2^13))=647
H^16(G,L(3.2^14))=1159
H^17(G,L(3.2^15))=2080\end{verbatim}}
\end{minipage}
\hspace{0.5cm}
\begin{minipage}[b]{0.5\linewidth}{\footnotesize
\centering
\begin{verbatim}
H^18(G,L(3.2^16))=3730
H^19(G,L(3.2^17))=6689
H^20(G,L(3.2^18))=12001
H^21(G,L(3.2^19))=21528
H^22(G,L(3.2^20))=38619
H^23(G,L(3.2^21))=69287
H^24(G,L(3.2^22))=124304
H^25(G,L(3.2^23))=223010
H^26(G,L(3.2^24))=400108
H^27(G,L(3.2^25))=717838
H^28(G,L(3.2^26))=1287890
H^29(G,L(3.2^27))=2310651
H^30(G,L(3.2^28))=4145619
H^31(G,L(3.2^29))=7437818
H^32(G,L(3.2^30))=13344508\end{verbatim}}
\end{minipage}
\end{figure}
The combinatorics become more complicated when one changes the value of $r$ away from $1$, though proofs of exponentiality using the above methods are available. One notices from the numbers, though, that the dimensions appear to grow at about the same rate as $1.8^m \sim 3.2^{m/2}$. 
\end{remark}

\begin{remark}\label{changep} For $p>2$ one can use essentially the same method to show that the sequence $\{\dim H^m(SL_2,L_m)\}$ also has exponential growth, where $L_m=L(2.p^m)$. 

We outline the changes necessary to show this:

The relevant recursions are \begin{align}\label{even}\Ext^q_G(\Delta(pb+i),M^{[1]}\otimes L(i))&\cong\bigoplus_{n\text{ even},\ 0\leq n\leq q}\Ext^{q-n}(\Delta(b+n),M)\\
\Ext^q_G(\Delta(pb+i),M^{[1]}\otimes L(\bar i))&\cong\bigoplus_{n\text{ odd},\ 0\leq n\leq q}\Ext^{q-n}(\Delta(b+n),M)\\
\Ext^q_G(\Delta(pb+p-1),M^{[1]}&\otimes L(p-1))\cong\Ext^q_G(\Delta(b),M)\end{align}
where $0\leq i\leq p-2$ and $\bar i=p-2-i$.

We use just equation (\ref{even}) above, starting with $b=i=0$. Then one continues to expand terms of the form $\Ext^{q}(\Delta(s),L_m)$ provided $p|s$ and $q$ is even and counts $\Ext^0$-leaves as before. 

Take in fact $m=2m'$; then an appropriate $a$-string $(a_1,\dots,a_m)$ with $\sum a_i=m$  is one for which \[\left(\frac{\frac{\frac{\frac{a_1}{p}+a_2}{p}+\dots}{\ddots}\ddots+a_{r-1}}{p}+a_r\right)\] is an integer for each $r\leq m$, where every $a_i$ is even and $2.p^m=\sum a_ip^i$. The continued fraction's integrality condition is equivalent to finding a $b$-string subject to $a_1=3b_1$ and $a_i+b_{i-1}=pb_i$ for each $i<m$; this also implies that each $b_i$ with $i<m$ is even. Interpreting the other restraints, we see such a $b$-string also satisfies $pb_i\geq b_{i-1}$ for $2\leq i\leq n-1$ and set $b_m=a_m$. We want that $m=\sum a_i=(p-1)\sum_{i=1}^{m-1} b_i+b_{n-1}+b_n$ with also $b_{n-1}+b_n=2$. Any string of non-negative integers satisfying these properties will work to give an $a$-string. One can then cook up exponentially many $b$-strings in a similar way to that done for $p=2$.
\end{remark}

\begin{remark} We have used Parker's equations to show that there is a sequence of simple modules $L_m$ with the value of $\dim \Ext^m_G(\Delta(0),L_m)$ growing exponentially. One can show similarly that there is a sequence $M_m$ with $\dim \Ext^m_G(\Delta(r),M_m)$ growing exponentially for any $r$. In fact, if $r<p^s$ then it is easy to see that $M_m=L_m^{[s]}\otimes L(r)$ will work. (One uses the fact that $\Ext^{q}_G(\Delta(pb+i),L(i)\otimes M^{[1]})\geq \Ext^q_G(\Delta(b),M)$.)\end{remark}

\begin{remark}\label{changeg}Brian Parshall asked by private communication if one could get exponential sequences $\{\dim H^m(G,L_m)\}$ for other $G$. We believe the answer is probably `yes' but as yet cannot give such a sequence. However we make some hopefully promising observations:

Firstly, let $G$ be any simple algebraic group with torus $T$. If $\lambda,\ \mu\in X^+(T)$ with $\lambda-\mu=m\beta$ for some $m\in \mathbb Z$ and $\beta$ a simple root, then, as observed in \cite{Par07}, we have by \cite[Corollary 10]{CPS04},

$$\label{eq}\Ext_G^q(\Delta(\lambda),L(\mu))\cong \Ext^q_{SL_2}(\Delta(2m_\beta),L(2n_\beta)),$$

where $m_\beta=\langle\lambda,\beta\rangle$ and $n_\beta=\langle \mu,\beta\rangle$. 

Now take $G=SL_3$ and $p=2$. We choose $\lambda_m=(2^m,0)$ on the $\alpha$-wall of the dominant chamber, where $\alpha=(2,-1)$ and $\beta=(-1,2)$ are the simple roots for $SL_3$ as elements of $X(T)$. Then take $\mu_m=\lambda_m+2^m\beta=(0,2^{m+1})$ and observe that in (\ref{eq}), with $\lambda=\lambda_m$, and $\mu=\mu_m$ we have $m_\beta=0$  with $n_\beta=2^{m+1}$.

Then we know from Remark \ref{changeweight} that the right hand side of \ref{eq} grows exponentially. 

If one knew that the number of composition factors $M_m$ of $\Delta(\lambda_m)$ admitting non-zero values of $\Ext^m(M_m,L(\mu_m))$ grew subexponentially, then one could find a sequence of such $M_m$ with the dimension of this latter $\Ext$ group growing exponentially. Since $M_m^*\otimes L(\mu_m)$ is irreducible by Steinberg, we would then have $\dim H^m(G,M_m^*\otimes L(\mu_m))$ giving the desired result. Unfortunately, using \cite[Theorem 4.12]{Par01} one can show there are $2^{m-2}+2$ composition factors in $\Delta(\lambda_m)$.\end{remark}

\begin{remark} While the example above doesn't give the exponential growth of cohomology asked for by Parshall, the same equation shows that for all $G$ and all $p$ we can take $\lambda$ and $\mu$ such that $\Ext^q_{SL_2}(\Delta(2m_\beta),L(2n_\beta))$ is big. This at least gives us that $\dim\Ext^q_G(\Delta(\lambda),L(\mu))$ has exponential behaviour as $\lambda$ and $\mu$ vary over all weights of $G$.\end{remark}

\begin{remark} It is remarkable that the dimensions of the modules in our sequences $\{L_m\}$ for which we have exponential growth of $H^m(G,L_m)$ are so small: when $G=SL_2$ and $p=2$, in Theorem \ref{mine} we used Frobenius twists of the two-dimensional natural module. Similarly, we could use three-dimensional modules when $p>2$.\end{remark}

This brings to mind some of the questions raised in \cite{GKKL}. We list some apposite results from that paper:
\begin{ttt}\begin{enumerate}\item Let $G$ be a finite simple group, $F$ a field and $M$ an $FG$ module. Then $\dim H^2(G,M)\leq 17.5\dim M$.
\item \label{b}Let $G$ be a finite group, $F$ a field and $M$ an irreducible $FG$ module. Then $\dim H^2(G,M)\leq 18.5\dim M$.
\item Let $F$ be an algebraically closed field of characteristic $p>0$ and $k$ a positive integer. Then there exists a sequence of finite groups $G_i, i\in \mathbb N$ and irreducible faithful $FG_i$-modules $M_i$ such that \begin{enumerate}\item $\lim_{i\to\infty}\dim M_i=\infty$,
\item $\dim H^k(G_i,M_i)\geq e(\dim M_i)^{k-1}$ for some constant $e=e(k,p)>0$, and
\item if $k\geq 3$ then $\lim_{i\to\infty}\frac{\dim H^k(G_i,M_i)}{\dim M_i}=\infty$.\end{enumerate}\end{enumerate}\end{ttt}
It is then pointed out that (iii) above precludes the possibility of generalising item (ii) above to higher degrees of cohomology. Nonetheless,  following questions are raised.

\begin{questions}\begin{enumerate}\item For which $k$ is it true that there is an absolute constant $C_k$ such that $\dim H^k(G,V) < C_k$ for all absolutely irreducible $FG$-modules $V$ and all finite 
simple groups $G$ with $F$ an algebraically closed field (of any characteristic)?
\item For which positive integers $k$ is it true that there is an absolute constant $d_k$ such that $\dim H^k(G,V) < d_k.\dim V^{k-1}$ for all absolutely irreducible faithful $FG$-modules $V$ and all finite groups $G$ with $F$ an algebraically closed field (of any characteristic)?\end{enumerate}\end{questions}

Note that there is no answer to question (i), for any $k>0$, even in the possibly easier case where $G$ is a simple algebraic group. The highest value of $\dim H^1(G,V)$ on record (see \cite{Sco03}) is $3$, where $G=SL_6$. Assuming Lusztig's character formula holds, we could take $p=7$ and $V=L(45454)$ to achieve this value. If we did have a positive answer to Question (i), this would imply a positive answer to Question (ii) in the case $G$ is taken to be a finite simple group.

In any case, our examples are relevant to Question (ii), when $G$ is taken to be a simple group. Consider the case when $G$ is algebraic. If $G$ is $SL_2$ we believe that $\max_{p,L-\text{irred}} \dim H^m(G,L)\leq \Pi_{m-1}$ with equality occurring if and only if $p=2$ and $L$ is a sufficiently high twist of $L(1)$. Then for all $G$, it is conceivable, owing to the low dimensions of the module involved, that the largest value of $\dim H^k(G,V)/(\dim V)^{k-1}$ occurs in the case $G=SL_2$, $p=2$ and where $V=L(1)^{[r]}$ is a twist of the natural module for $G$, since then, $\dim V=2$ and the lowest it could possibly be. But while the rate of growth of $\max_r \dim H^k(G,L(1)^{[r]})$ is exponential, it grows at about the rate $1.8^k$, so that $\dim H^k(G,V)/(\dim V)^{k-1}\sim 1.8^k/2^{k-1}$ will tend to zero. 

Thus it is conceivable that one could ask for a single constant $d\geq d_k$ that works for all $k$ in Question (ii), when $G$ is simple and algebraic. Ignoring the case where $k=1$ (and Questions (i) and (ii) coincide), possibly even $d=1$ may work. This is then relevant to the finite group situation by considering generic cohomology. One has from \cite{CPSK} that $H^m(G,V^{[e]})\cong H^m(G(q),V)$ for high enough values of $e$ and $q$. Our example provides some small evidence then, that for $k>1$, one might replace $d_k$ with a universal constant in Question (ii) if $G$ is a finite simple group.

{\footnotesize\bibliographystyle{amsalpha}
\bibliography{psans.bib}
}
\end{document}